\documentclass[10pt]{amsart}
\usepackage{natbib}
\usepackage{amssymb}
\theoremstyle{definition}
\newtheorem{definition}{Definition}

\theoremstyle{plain}
\newtheorem{theorem}[definition]{Theorem}
\newtheorem{lemma}[definition]{Lemma}

\def\oo{\sigma}
\def\pp{\tau}
\def\ii{\iota}
\def\bb#1{\overline{#1}}
\def\nn#1{#1^{-1}}
\def\xx{x^{-1}}
\def\yy{y^{-1}}
\def\mm#1#2#3#4{\left(\begin{array}{cc}#1&#2\\#3&#4\end{array}\right)}
\def\opos#1{#1^{\textrm{op}}}
\def\a{\alpha}
\def\b{\beta}
\def\c{\gamma}
\def\d{\delta}
\def\e{\varepsilon}

\title[Loops $M(G,2)$]
{On the Uniqueness of Loops $M(G,2)$}

\author[Petr~Vojt\v{e}chovsk\'y]{Petr Vojt\v{e}chovsk\'y}
\email{petr@math.du.edu}
\address{Department of Mathematics, University of Denver, 2360 S Gaylord St,
Denver, CO 80208, USA}

\thanks{Work partially supported by Grant Agency of Charles University, grant number
$269/2001/$B-MAT/MFF}

\begin{document}

\begin{abstract}Let $G$ be a finite group and $C_2$ the cyclic group of order
$2$. Consider the $8$ multiplicative operations $(x,y)\mapsto (x^iy^j)^k$,
where $i$, $j$, $k\in\{-1,\,1\}$. Define a new multiplication on $G\times
C_2$ by assigning one of the above $8$ multiplications to each quarter
$(G\times\{i\})\times(G\times\{j\})$, for $i$, $j\in C_2$. When $G$ is
nonabelian then exactly four assignments yield Moufang loops that are not
associative; all $($anti$)$isomorphic, known as loops $M(G,2)$. \vskip 1mm

 \vskip 1mm

\end{abstract}

\maketitle

\section{Introduction}

\noindent Because of the specialized topic of this paper, we assume that the
reader is familiar with the theory of Bol and Moufang loops (cf.\
\cite{Pflugfelder}).

Chein introduced the following construction in \cite{Chein} to obtain Moufang
loops from groups: Let $G$ be a finite group and let $\bb{G}=\{\bb{x};\; x\in
G\}$ be a set of new elements. Define multiplication $*$ on $G\cup\bb{G}$ by
\begin{equation}\label{Eq:Chein}
    x*y = xy,\quad x*\bb{y}=\bb{yx},\quad
    \bb{x}*y=\bb{x\yy},\quad \bb{x}*\bb{y}=\yy x,
\end{equation}
where $x$, $y\in G$. The resulting Moufang loop $M(G,2)$ is associative if and
only if $G$ is abelian, according to \cite{Chein}.

Loops $M(G,2)$ play an important role among Moufang loops of small order (cf.\
\cite{Chein}, \cite{GMR}). Recently, it was found that all Moufang loops of
order $n\in\{16,\,24,\,32\}$ can be obtained by modifying one quarter of the
multiplication tables of loops $M(G,\,2)$ in a certain way \cite{AlesPetr}. The
smallest nonassociative Moufang loop is isomorphic to $M(S_3,2)$, where $S_3$
is the symmetric group on $3$ points (cf.\ \cite{CP}, \cite{Vojtechovsky}).

We are going to study a generalization of Chein's construction
(\ref{Eq:Chein}). Given a group $G$, consider the $8$ multiplicative operations
on $G$: $(x,y)\mapsto (x^iy^j)^k$, where $i$, $j$, $k\in \{-1,1\}$. Let $C_2$
be the cyclic group of order $2$. Define a new multiplication on $G\times C_2$
by assigning one of the above $8$ multiplications to each quarter
$(G\times\{i\})\times(G\times\{j\})$, for $i$, $j\in C_2$. Let $M$ be the
resulting quasigroup.

It this note, we characterize when $M$ is a loop (Lemma \ref{Lm:One}); we show
that if $M$ is a Bol loop, it is Moufang (Lemma \ref{Lm:IP}); and we prove that
when $G$ is nonabelian then there are exactly $4$ assignments that yield
nonassociative Moufang loops, all (anti)isomorphic to the loop $M(G,2)$. See
Theorem $6$ for details.

Chein's construction $(\ref{Eq:Chein})$ is therefore unique, in a sense.

\section{Notation}

\noindent Let us introduce a notation that will better serve our purposes.
Consider the permutations $\ii$, $\oo$, $\pp$ of $G\times G$ defined by
$(x,y)\ii=(x,y)$, $(x,y)\oo=(y,x)$, and $(x,y)\pp=(\yy,x)$. Since
$\oo^2=\pp^4=\ii$ and $\oo\pp\oo=\nn{\pp}$, the group $A$ generated by $\oo$
and $\pp$ is isomorphic to $D_8$, the dihedral group of order $8$. The
elements $\psi$ of $A$ are described by
\begin{footnotesize}\begin{displaymath}
\begin{array}{c|cccccccc}
    \psi&\ii&\oo&\pp&\pp^2&\pp^3&\oo\pp&\oo\pp^2&\oo\pp^3\\
    (x,y)\psi& (x,y)&(y,x)&(\yy,x)&(\xx,\yy)&(y,\xx)&(\xx,y)&
        (\yy,\xx)&(x,\yy).
\end{array}
\end{displaymath}\end{footnotesize}\noindent
We like to think of these elements as multiplications in $G$, and often
identify $\psi\in A$ with the map $\psi\Delta:G\times G\to G$, where
$(x,y)\Delta=xy$. For instance, the permutation $\oo\pp$ determines the
multiplication $x*y=\xx y$. Note that $\oo\Delta=\ii\Delta$ when $G$ is
abelian, and that $A\Delta=\ii\Delta$ when $G$ is an elementary abelian
$2$-group.

To avoid trivialities, we assume throughout the paper that \emph{$G$ is not an
elementary abelian $2$-group, and that $|G|>1$}.

It is natural to split the multiplication table of $M(G,2)$ into four quarters
$G\times G$, $G\times\bb{G}$, $\bb{G}\times G$ and $\bb{G}\times\bb{G}$, as in
\begin{displaymath}
\begin{array}{c|cc}
    *&G&\bb{G}\\
    \hline
    G&&\\
    \bb{G}&&
\end{array}\;.
\end{displaymath}
Then Chein's construction (\ref{Eq:Chein}) can be represented by the matrix
\begin{equation}\label{Eq:MatrixChein}
    M_c=\mm{\ii}{\oo}{\oo\pp^3}{\pp}.
\end{equation}
For example, we can see from $M_c$ that $\bb{x}*y=\bb{(x,y)\oo\pp^3} =
\bb{x\yy}$, for $x$, $y\in G$.

\section{Main Result}

\noindent When we look at Chein's construction $(\ref{Eq:Chein})$ via
$(\ref{Eq:MatrixChein})$, it appears to be somewhat arbitrary. Let us therefore
investigate all multiplications
\begin{equation}\label{Eq:GenericLoop}
    M=\mm{\a}{\b}{\c}{\d},
\end{equation}
where $\a$, $\b$, $\c$, $\d\in A$. We will no more distinguish between the
matrix $M$ and the groupoid it defines.

We note in passing that every $M$ is a quasigroup. The next Lemma characterizes
all loops $M$. In the course of the proof we encounter several identities of
the form $w_1=w_2$, where $w_i$ is a word in some symbols $x_1$, $\dots$,
$x_m\in G$. When $w_1$, $w_2$ reduce to the same word in the free group on
$x_1$, $\dots$, $x_m$, then $w_1=w_2$ surely holds in $G$. Conversely, since we
assumed that $G$ is not an elementary abelian $2$-group and $|G|>1$, there are
many identities that do not hold in $G$, no matter what $G$ is. For instance,
$x\ne x^{-1}$, $y\ne xy^{-1}x^{-1}$ (set $x=y$), and so on.

\begin{lemma}\label{Lm:One}
$M$ is a loop if and only if $\a\in\{\ii$, $\oo\}$, $\b\in\{\ii$, $\oo$,
$\pp^3$, $\oo\pp\}$ and $\c\in\{\ii$, $\oo$, $\pp$, $\oo\pp^3\}$. When $M$ is a
loop, its neutral element coincides with the neutral element of $G$.
\end{lemma}
\begin{proof}
We first show that if $M$ is a loop, its neutral element $e$ coincides with
the neutral element $1$ of $G$. This is clear, as for some $\e\in A$ we have
$1=1*1=(1,\,1)\e = 1 = 1*e$, and thus $1=e$.

The equation $y=1*y$ holds for every $y\in G$ if and only if $y=(1,y)\alpha$,
which happens if and only if $\a\in\{\ii$, $\oo$, $\pp^3$, $\oo\pp\}$.
Similarly, the equation $y=y*1$ holds for every $y\in G$ if and only if
$\a\in\{\ii$, $\oo$, $\pp$, $\oo\pp^3\}$. Altogether, $y=y*1=1*y$ holds for
every $y\in G$ if and only if $\a\in\{\ii$, $\oo\}$.

Following the same strategy, $\bb{y}=1*\bb{y}$ holds for every $y\in G$ if and
only if $\b\in\{\ii$, $\oo$, $\pp^3$, $\oo\pp\}$, and $\bb{y}=\bb{y}*1$ holds
for every $y\in G$ if and only if $\c\in\{\ii$, $\oo$, $\pp$, $\oo\pp^3\}$.
\end{proof}

Once $M$ is a loop, it must have two-sided inverses:

REVISION: IN THE ORIGINAL VERSION I CLAIMED THAT IF $M$ IS A LOOP THAT IS BOL
THEN IT MUST BE MOUFANG. IT'S NOT TRUE.

\begin{lemma}\label{Lm:IP}
If $M$ is a loop then it is an inverse property loop.
\end{lemma}
\begin{proof}
Assume that $x*y=1$ for some $x$, $y\in G\cup \bb{G}$. Then both $x$, $y$
belong to $G$, or both belong to $\bb{G}$, by Lemma \ref{Lm:One}. We therefore
want to show that $(x,y)\e=1$ implies $(y,x)\e=1$ for every $\e\in A$ and $x$,
$y\in G$.

Pick $\e\in A$. Then $(x,y)\e=(x^iy^j)^k$ for some $i$, $j$, $k\in\{-1,\ 1\}$.
Assume that $(x,y)\e=1$. Then $x^iy^j=1$ and $y^jx^i=1$. If $i=j$, we conclude
from the latter equality that $y^ix^j=1$, and thus $(y,x)\e=1$. The inverse of
the former equality yields $y^{-j}x^{-i}=1$. If $i=-j$, we immediately have
$y^ix^j=1$, and thus $(y,x)\e=1$.

Hence $M$ is an inverse property loop.
\end{proof}

Given $M$ as in $(\ref{Eq:GenericLoop})$, let
\begin{displaymath}
    \opos{M}=\mm{\oo\a}{\oo\c}{\oo\b}{\oo\d}.
\end{displaymath}

\begin{lemma}\label{Lm:w}
The quasigroup $\opos{M}$ is opposite to $M$.
\end{lemma}
\begin{proof}
Denote by $\circ$ the multiplication in $\opos{M}$. Then
\begin{eqnarray*}
    x\circ y&=&(x,y)\oo\a=(y,x)\a=y*x,\\
    x\circ\bb{y}&=&\bb{(x,y)\oo\c}=\bb{(y,x)\c}=\bb{y}*x,\\
    \bb{x}\circ y&=&\bb{(x,y)\oo\b}=\bb{(y,x)\b}=y*\bb{x},\\
    \bb{x}\circ\bb{y}&=&(x,y)\oo\d=(y,x)\d=\bb{y}*\bb{x},
\end{eqnarray*}
for every $x$, $y\in G$.
\end{proof}

Let us assume from now on that \emph{$G$ is nonabelian}. Then the identity
$xy=yx$ and any other identity that reduces to $xy=yx$ do not hold in $G$, of
course. We will come across the identity $xxy=yxx$. Note that this identity
holds in $G$ if and only if the center of $G$ is of index $2$ in $G$.

We would like to know when $M$ is a Bol (and hence Moufang) loop. Assume from
now on that $M$ is a loop.

Recall that the opposite of a Moufang loop is again Moufang. We can therefore
combine Lemmas \ref{Lm:One}, \ref{Lm:w} and assume that the loop $M$ satisfies
$\a=\ii$. Since every Moufang loop is diassociative, we are going to have a
look at such loops:

\begin{lemma}\label{Lm:Diass}
If $G$ is nonabelian and $M$ is a diassociative loop with $\a=\ii$ then
$(\b,\c,\d)$ is one of the eight triples
\begin{equation}\label{Eq:Diass2}
    \begin{array}{llll}
        (\ii,\ii,\ii),&(\pp^3,\ii,\oo\pp),&
        (\oo,\oo,\oo),& (\oo\pp,\oo,\pp^3),\\
        (\pp^3,\pp,\pp^2),& (\ii,\pp,\oo\pp^3),&
        (\oo,\oo\pp^3,\pp),& (\oo\pp,\oo\pp^3,\oo\pp^2).
    \end{array}
\end{equation}
\end{lemma}
\begin{proof}
The identities $(\bb{x}*\bb{x})*y=\bb{x}*(\bb{x}*y)$,
$\bb{x}*(y*\bb{x})=(\bb{x}*y)*\bb{x}$ hold in $M$, for every $x$, $y\in G$.
They translate into
\begin{eqnarray}
    (x,x)\d y&=& (x,(x,y)\c)\d\label{Eq:Diass1},\\
    (x,(y,x)\b)\d&=&((x,y)\c,x)\d\label{Eq:Flex},
\end{eqnarray}
respectively. We are first going to check which pairs $(\c,\d)$ satisfy
$(\ref{Eq:Diass1})$.

Assume that $\c=\ii$. Then $(\ref{Eq:Diass1})$ becomes $(x,x)\d y=(x,xy)\d$.
Denote this identity by $I(\delta)$. Then $I(\ii)$ is $xxy=xxy$ (true),
$I(\oo)$ is $xxy=xyx$ (false), $I(\pp)$ is $y=y^{-1}$ (false), $I(\pp^2)$ is
$x^{-2}y=x^{-1}y^{-1}x^{-1}$ (false), $I(\pp^3)$ is $y=xyx^{-1}$ (false),
$I(\oo\pp)$ is $y=y$ (true), $I(\oo\pp^2)$ is $x^{-2}y=y^{-1}x^{-1}x^{-1}$
(false), and $I(\oo\pp^3)$ is $y=xy^{-1}x^{-1}$ (false).

Assume that $\c=\oo$. Then $(\ref{Eq:Diass1})$ becomes $(x,x)\d y=(x,yx)\d$.
Verify that this identity holds only if $\delta=\oo$ or $\delta=\pp^3$. (The
case $\delta=\oo$ leads to the identity $xxy=yxx$ mentioned before this Lemma.)

When $\c=\pp$, $(\ref{Eq:Diass1})$ holds only if $\d=\pp^2$ or $\d=\oo\pp^3$.

When $\c=\oo\pp^3$, $(\ref{Eq:Diass1})$ holds only if $\d=\pp$ or
$\d=\oo\pp^2$.

Altogether, $(\ref{Eq:Diass1})$ can be satisfied only when $(\c,\d)$ is one of
the $8$ pairs $(\ii,\ii)$, $(\ii,\oo\pp)$, $(\oo,\oo)$, $(\oo,\pp^3)$,
$(\pp,\pp^2)$, $(\pp,\oo\pp^3)$, $(\oo\pp^3,\pp)$, $(\oo\pp^3,\oo\pp^2)$. All
these pairs will now be tested on $(\ref{Eq:Flex})$.

Straightforward calculation shows that $(\ref{Eq:Flex})$ can be satisfied only
when $(\b,\c,\d)$ is one of the $8$ triples listed in $(\ref{Eq:Diass2})$.
\end{proof}

The \emph{Moufang identity} $((xy)x)z=x(y(xz))$ will help us eliminate $4$ out
of the $8$ possibilities in $(\ref{Eq:Diass2})$. We have
$((x*\bb{y})*x)*z=x*(\bb{y}*(x*z))$ in $M$, and thus
\begin{equation}\label{Eq:4}
    (((x,y)\b,x)\c,z)\c = (x,(y,xz)\c)\b.
\end{equation}
The pairs $(\b,\c)=(\oo,\oo)$, $(\pp^3,\ii)$, $(\ii,\pp)$, $(\oo\pp,\oo\pp^3)$
do not satisfy $(\ref{Eq:4})$. For instance, $(\b,\c)=(\oo,\oo)$ turns
$(\ref{Eq:4})$ into $zxyx=xzyx$, i.e., $zx=xz$.

The four remaining triples from $(\ref{Eq:Diass2})$ yield Moufang loops, as we
are going to show.

The quadruple $(\alpha,\beta,\gamma,\delta) = (\ii,\ii,\ii,\ii) =
G_\ii$ corresponds to the direct product of $G$ and the two-element
cyclic group. The quadruple $(\ii,\oo,\oo\pp^3,\pp)=M_c$ is the
Chein Moufang loop $M(G,2)$ that is associative if and only if $G$
is abelian, by \cite{Chein}. (We can also verify this directly.)

Set $G_\pp=(\ii,\pp^3,\pp,\pp^2)$ and $M_\oo=(\ii,\oo\pp,\oo,\pp^3)$. We claim
that $G_\ii$ is isomorphic to $G_\pp$, and $M_c$ is isomorphic to $M_\oo$.

\begin{lemma}\label{Lm:T} Define $T:A^4\to A^4$ by
\begin{displaymath}
    M=\mm{\a}{\b}{\c}{\d}\mapsto\mm{\a}{\pp^3\b}{\c\pp}{\pp^2\d}=MT.
\end{displaymath}
If $((x,y)\b\Delta)^{-1}=(y^{-1},x^{-1})\b\Delta$ and
$((x,y)\c\Delta)^{-1}=(x^{-1},y)\c\pp\Delta$ then $M$ is isomorphic to $MT$.
\end{lemma}
\begin{proof}
Consider the permutation $f$ of $G\cup \bb{G}$ defined by $f(x)=x$,
$f(\bb{x})=\bb{x^{-1}}$, for $x\in G$. Let $*$ be the multiplication in $M$ and
$\circ$ the multiplication in $MT$. We show that $(x*y)f=xf\circ yf$ for every
$x$, $y\in G\cup\bb{G}$. With $x$, $y\in G$, we have
\begin{eqnarray*}
    &&(x*y)f=(x,y)\a\Delta f=(x,y)\a\Delta=x\circ y=xf\circ yf,\\
    &&(\bb{x}*\bb{y})f = (x,y)\d\Delta f= (x,y)\d\Delta
    =(x^{-1},y^{-1})\pp^2\d\Delta = \bb{x}f\circ \bb{y}f.
\end{eqnarray*}
Using the assumption on $\beta$ and $\gamma$, we also have
\begin{displaymath}
    (x*\bb{y})f=\bb{(x,y)\b\Delta}f=\bb{((x,y)\b\Delta)^{-1}}
    =\bb{(y^{-1},x^{-1})\b\Delta}=\bb{(x,y^{-1})\pp^3\b\Delta} = xf\circ \bb{y}f,
\end{displaymath}
and
\begin{displaymath}
    (\bb{x}*y)f=\bb{(x,y)\c\Delta}f=\bb{((x,y)\c\Delta)^{-1}}
    =\bb{(x^{-1},y)\c\pp\Delta}=\bb{x}f\circ yf.
\end{displaymath}
\end{proof}

Note that $G_\ii T=G_\pp$ and $M_cT=M_\oo$. Now, $\b\in\{\ii,\oo\}$ satisfies
$((x,y)\b\Delta)^{-1}=(y^{-1},x^{-1})\b\Delta$, and $\c\in\{\ii,\oo\pp^3\}$
satisfies $((x,y)\c\Delta)^{-1}=(x^{-1},y)\c\pp\Delta$. By Lemma \ref{Lm:T},
$G_\ii$ is isomorphic to $G_\pp$, and $M_c$ is isomorphic to $M_\oo$.

We have proved:

\begin{theorem}\label{Th:Main}
Let $G$ with $|G|>1$ be a finite group that is not an elementary abelian
$2$-group. With the above conventions, let
\begin{displaymath}
    M=\mm{\a}{\b}{\c}{\d}
\end{displaymath}
specify the multiplication in $L=G\cup\bb{G}$, where $\a$, $\b$, $\c$, $\d\in
A=\langle\oo,\pp\rangle$, and $(x,y)\oo=(y,x)$, $(x,y)\pp=(y^{-1},x)$.

When $G$ is nonabelian, then $L$ is a Moufang loop if and only if $M$ is
equal to one of the following matrices:
\begin{displaymath}
    \begin{array}{ll}
    G_\ii=\mm{\ii}{\ii}{\ii}{\ii},&
    \opos{G_\ii}=\mm{\oo}{\oo}{\oo}{\oo},\\
    G_\pp=\mm{\ii}{\pp^3}{\pp}{\pp^2},&
    \opos{G_\pp}=\mm{\oo}{\oo\pp}{\oo\pp^3}{\oo\pp^2},\\
    M_c=\mm{\ii}{\oo}{\oo\pp^3}{\pp},&
    \opos{M_c}=\mm{\oo}{\pp^3}{\ii}{\oo\pp},\\
    M_\oo=\mm{\ii}{\oo\pp}{\oo}{\pp^3},&
    \opos{M_\oo}=\mm{\oo}{\ii}{\pp}{\oo\pp^3}.
    \end{array}
\end{displaymath}
The loops $\opos{X}$ are opposite to the loops $X$. The isomorphic loops
$G_\ii$, $G_\pp$ and their opposites are groups. The isomorphic loops $M_c$,
$M_\oo$ and their opposites are Moufang loops that are not associative.
\end{theorem}

\bibliographystyle{plain}

\begin{thebibliography}{99}

\bibitem[Chein(1978)]{Chein} Orin~Chein, \emph{Moufang loops of small order},
Memoirs of the American Mathematical Society, Volume \textbf{13}, Issue 1,
Number \textbf{197} (1978).

\bibitem[Chein, Pflugfelder (1971)]{CP} O.~Chein, H.~O.~Pflugfelder, \emph{The
smallest Moufang loop}, Arch.\ Math.\ \textbf{22} (1971), 573--576.

\bibitem[Drapal and Vojtechovsky(2002)]{AlesPetr} Ale\v{s} Dr\'apal and Petr
Vojt\v{e}chovsk\'y, \emph{Moufang loops that share associator and three
quarters of their multiplication tables}, submitted.

\bibitem[GMR(1999)]{GMR} Edgar G.~Goodaire, Sean May, Maitreyi Raman, \emph{The
Moufang Loops of Order less than $64$}, Nova Science Publishers, 1999.

\bibitem[Pflugfelder(1990)]{Pflugfelder}
H.~O.~Pflugfelder, Quasigroups and Loops: Introduction, \emph{Sigma series in
pure mathematics} {\bf 7}, Heldermann Verlag Berlin, 1990.

\bibitem[Vojtechovsky(2003)]{Vojtechovsky} Petr Vojt\v{e}chovsk\'y, \emph{The
smallest Moufang loop revisited}, to appear in Results in
Mathematics.

\end{thebibliography}

\end{document}